\documentclass[11pt]{article}
\usepackage[T1]{fontenc}
\usepackage{lmodern}
\usepackage{amsmath,amssymb,amsthm}
\usepackage{comment}
\usepackage[a4paper,margin=29mm]{geometry}
\usepackage{microtype}
\usepackage[hidelinks]{hyperref}

\newtheorem{theorem}{Theorem}
\newtheorem{lemma}[theorem]{Lemma}
\newtheorem{proposition}[theorem]{Proposition}
\newtheorem{corollary}[theorem]{Corollary}
\newtheorem{claim}{Claim}[theorem]
\newcommand{\len}{\ell}
\newcommand{\J}{\{0,1,2,3\}}

\title{An improved bound for the intersection of two longest cycles in \(k\)-connected graphs}
\author{Juan Guti\'errez}

\begin{document}
\maketitle

\begin{abstract}
A conjecture attributed to Smith states that in every $k$-connected
graph, every pair of longest cycles intersects in at least $k$ vertices.
We prove that every two longest cycles in a $k$-connected graph on $n$
vertices intersect in at least $\min\{k,10k-n-30\}$ vertices. Consequently,
Smith's conjecture holds whenever $k\ge \frac{n+30}{9}$. This improves the previous bound of $\frac{n+16}{7}$ 
for $n>33$ and also implies a
bound of $\frac{n+22}{9}$ for paths instead of cycles.
\end{abstract}

\section{Introduction}

Every two longest paths in a connected graph share a vertex, and every two
longest cycles in a $2$-connected graph share at least two vertices. A
natural question is whether higher connectivity forces larger
intersections. A conjecture attributed to Scott Smith and recorded by
Gr\"otschel~\cite[Conjecture~5.2]{Gro84} asserts that, for every $k\ge2$,
any two longest cycles in a $k$-connected graph intersect in at least $k$
vertices. Hippchen~\cite[Conjecture~2.2.4]{Hippchen} later proposed the
corresponding statement for longest paths. For background on intersections
of longest paths and cycles, we refer to the survey of Shabbir,
Zamfirescu, and Zamfirescu~\cite{SZZ}.

Smith's conjecture is known for $2\le k\le8$, by results of Gr\"otschel
and Nemhauser~\cite{GN84}, Gr\"otschel~\cite{Gro84}, and Stewart and
Thompson~\cite{ST95}. For arbitrary $k$, Chen, Faudree, and
Gould~\cite{CFG} proved that every two longest cycles intersect in at
least $c_0k^{3/5}$ vertices, where
$c_0=(\sqrt[3]{256}+3)^{-3/5}\approx0.2615$.
Recent work has substantially improved this general lower bound.
Groenland, Longbrake, Steiner, Turcotte, and Yepremyan~\cite{GLSTY}
obtained $\Omega(k^{5/8})$, and Ma and Zhao~\cite{MZ} improved it to
$\Omega(k^{2/3})$. In a subsequent preprint, D.~M.~Chen~\cite{DChen}
established $\Omega(k^{8/11})$. Most recently, Ma, Ning, and
Zhao~\cite{MNZ} proved that every two longest cycles share
at least $k/600$ vertices, giving the first general linear lower bound.

For Hippchen's conjecture, the cases $k\in \{2,3\}$ were proved by
Hippchen~\cite{Hippchen}, the case $k=4$ by the author~\cite{G2021}, and the case $k=5$ by Cho, Choi, and
Park~\cite{CCP}. Guti\'errez and Valqui~\cite{GV} subsequently verified
the conjecture for every $1\le k\le7$ by relating it to Smith's
conjecture. The recent linear cycle bound also yields a general lower
bound of $\lfloor k/600\rfloor$ for the intersection of two longest
paths~\cite[Theorem~9.2]{MNZ}. Related results for cycles and paths whose
lengths are close to the maximum, but which need not themselves be
longest, were obtained by Wu~\cite{Wu2026}.

 The author showed that 
 Hippchen’s conjecture is valid when 
 $k \geq \frac{n-2}{3}$ \cite{G2021}. 
Cho, Choi, and Park~\cite{CCP} improved this result, proving Hippchen's
conjecture whenever $k\ge \frac{n+2}{5}$. For cycles, Guti\'errez and
Valqui~\cite{GV} 
confirm Smith's conjecture whenever
$k\ge(n+16)/7$. Our main result improves this bound as follows.

\begin{theorem}\label{thm:main}
Let $G$ be a $k$-connected graph on $n$ vertices, where $k\ge2$, and let
$C,D$ be longest cycles of $G$. Then
$
 |V(C)\cap V(D)|\ge\min\{k,10k-n-30\}.
$
In particular, Smith's conjecture holds whenever $k\ge \frac{n+30}{9}$.
\end{theorem}

To state the corresponding result for paths, we recall the notation and
reduction introduced in~\cite{GV}. For integers $k\ge1$ and $1\le r\le k$,
let $P(k,r)$ be the statement that every two longest paths in a
$k$-connected graph intersect in at least $r$ vertices, and let $C(k,r)$
be the analogous statement for longest cycles. Thus $P(k,k)$ is
Hippchen's conjecture, while $C(k,k)$, for $k\ge2$, is Smith's conjecture.
Guti\'errez and Valqui~\cite[Observation~3]{GV} showed that
$C(k+1,r+1)$ implies $P(k,r)$.
Applying this to
Theorem~\ref{thm:main} gives the following corollary.

\begin{corollary}\label{cor:paths}
Let $G$ be a $k$-connected graph on $n$ vertices, $k\ge1$, and let $P,Q$
be longest paths of $G$. Then
$
 |V(P)\cap V(Q)|\ge\min\{k,10k-n-22\}.
$
In particular, Hippchen's conjecture holds whenever $k\ge \frac{n+22}{9}$.
\end{corollary}

The main ingredient in the proof of Theorem~\ref{thm:main} is a
four-fan argument. If a component outside a longest cycle contains a
path on four vertices, we prove that the cycle has length at least
$5(k-3)$. 
If neither exterior contains a path on four vertices, its
components are stars or triangles, and we obtain an
intersection of size at least $k$.

\begin{proof}[Proof of Theorem~\ref{thm:main}]
Put $L=|C|=|D|$ and $X=V(C)\cap V(D)$.
If an exterior component of either cycle contains a path on four vertices,
Lemma~\ref{lem:fourfan} gives $L\ge5(k-3)$. Hence
\[
 |X|=2L-|V(C)\cup V(D)|\ge2L-n\ge10k-n-30.
\]
Otherwise Lemma~\ref{lem:nop4} gives $|X|\ge k$.
The stated minimum follows. If $9k\ge n+30$, then $10k-n-30\ge k$,
which proves the final assertion.
\end{proof}

\section{Preliminaries}

All graphs are finite and simple, and all paths and cycles are simple. The
length $|P|$ of a path or cycle is its number of edges. For a fixed orientation of a cycle
$C$ and distinct vertices $u,v\in V(C)$, write $C[u,v]$ for the subpath from
$u$ to $v$ in the orientation. We write $G-C$ for $G-V(C)$ and
$N_C(H):=N_G(V(H))\cap V(C)$ for a component $H$ of $G-C$.

We use the following standard consequence of Menger's theorem, the Fan
Lemma:
\begin{proposition}[\cite{BM}]\label{prop:kfan}
    Let $G$ be a $k$-connected graph. 
    Let $S \subseteq V(G)$.
If $v\notin S$, and $|S|\ge k$, then there
are $k$ paths from $v$ to distinct vertices of $S$, each meeting
$S$ only at its last vertex, whose pairwise intersection is exactly
$\{v\}$.
\end{proposition}
We call such a collection a \emph{$k$-fan} \emph{from} $v$ \emph{to} $S$.



\begin{proposition}\label{prop:basica}
Let $C$ be a longest cycle and let $u,v$ be distinct vertices of $C$.
If $R$ is a $u$--$v$ path whose internal vertices lie outside $C$,
then $|C[u,v]|\ge |R|$.
\end{proposition}
\begin{proof}
Otherwise, replacing $C[u,v]$ by $R$ gives a longer cycle.
\end{proof}

\begin{proposition}\label{prop:attachments}
Let $C$ be a longest cycle of a $k$-connected graph, where $k\ge2$,
and let $H$ be a component of $G-C$. Then $|N_C(H)|\ge k$ and
$|C[u,v]|\ge2$ for all distinct $u,v\in N_C(H)$.
Consequently, if $C$ is not Hamiltonian, then $|C|\ge2k$.
\end{proposition}
\begin{proof}
Since $\delta(G)\ge k$, an end of a longest path, together with
its farthest neighbor on that path, gives a cycle of length at least
$k+1$. Thus $|C|\ge k+1$. If $|N_C(H)|<k$, then deleting $N_C(H)$
separates $H$ from a vertex of $C$, contradicting $k$-connectivity.
For any distinct $u,v\in N_C(H)$, connectedness of $H$ gives a
$u$--$v$ path of length at least two with all internal vertices in
$H$, so $|C[u,v]| \geq 2$ by Proposition \ref{prop:basica}.
For the last assertion, list the vertices of $N_C(H)$ as
$u_1,\ldots,u_p$ in cyclic order on $C$, and set $u_{p+1}=u_1$.
Summing $|C[u_i,u_{i+1}]|\ge2$ over $1\le i\le p$ gives
$|C|\ge2p\ge2k$.
\end{proof}



\begin{proposition}\label{prop:routing}
Let $C$ be a longest cycle with 2 distinct vertices $x$ and $y$.
Let $R$ be a path disjoint from $C$
with ends $a$ and $b$ (possibly $a=b$).
Let $P$ and $P'$ be two internally disjoint paths from
$a$ to $x$ and $y$, respectively.
Let $Q$ be a path from $b$ to $y$.
If $P,P'$  and $Q$ are internally disjoint from both $R$
and $C$, then $|C[x,y]| \ge 2+|R|$.
\end{proposition}
\begin{proof}
    Traverse $Q$ from $b$, and let $c$ be its first intersection with
$P\cup P'$. Such an intersection exists because $y\in P'$. 
If $c\in P$ then
we apply Proposition \ref{prop:basica}
with path $P[x,c] \cdot Q[c,b] \cdot R[b,a] \cdot P'$, obtaining $|C[x,y]|\geq 2+|R|$.
Otherwise, 
we apply Proposition \ref{prop:basica}
with path $P[x,a] \cdot R[a,b] \cdot Q[b,c] \cdot P'[c,y]$,
obtaining the same bound.
\end{proof}

\section{Proof of Lemma \ref{lem:fourfan}}

\begin{lemma}\label{lem:fourfan}
Let $C$ be a longest cycle of a $k$-connected graph. If a component $H$ of
$G-C$ contains a path $R=v_0v_1v_2v_3$, then
$
 |C|\ge5(k-3).
$
\end{lemma}
\begin{proof}
If $k\le5$ then, by Proposition \ref{prop:attachments}, $|C|\ge 2k \geq 5k-15$, so assume $k\ge6$.
For each $j\in\J$, take a $k$-fan from $v_j$ to $C$, which exists by Proposition \ref{prop:kfan}.
At most $3$ paths contain another vertex of $R$, since the paths in one
fan meet only at $v_j$. 
Retain, for each $j$, a
collection of $k-3$ paths
from $v_j$ to $C$,
which we call $\mathcal F_j$,
and let $\mathcal F = \mathcal F_0 \cup \mathcal F_1 \cup \mathcal F_2 \cup \mathcal F_3$.
Thus each path meets $R$ only at $v_j$ and $C$
only at its other end, and distinct paths in $\mathcal F_j$ meet
only at $v_j$. 
Also, every vertex of a path in $\mathcal F$ belongs to $H$ except for one of its ends, which is in $C$.

Moreover, letting $F$ be the union of the ends of the paths in $\mathcal F$ that are in $C$,
among all such $\mathcal F$ choose one minimizing $|F|$.
Define, for every $u \in F$,
\[
 S(u)=\{j:\text{a path of }\mathcal F_j\text{ ends at }u\}.
\]
Let $u_1,\ldots,u_f$ be the vertices of $F$ in the cyclic order they appear on $C$, and set
$ S_i=S(u_i)$. Let $\len_i=|C[u_i,u_{i+1}]|,
$
with indices modulo $f$.

\begin{claim}\label{cl:overlap}
If $S_i \cap S_{i+1} \neq \emptyset$, then
$
 \len_i\ge2+\max\{|a-b|:a\in S_i\cap S_{i+1} ,\ b\in S_i\cup S_{i+1}\}.
$
\end{claim}
\begin{proof}
Take
$a\in S_i\cap S_{i+1}$ and $b\in S_i\cup S_{i+1}$ attaining the maximum. The two
paths of $\mathcal F_a$ ending at $u_i,u_{i+1}$ meet only at $v_a$.
Take also the path of $\mathcal F_b$ ending at one of these two vertices,
and apply Proposition~\ref{prop:routing}, interchanging the ends if necessary.
\end{proof}


\begin{claim}\label{cl:shortbad}
If $S_i \cap S_{i+1} = \emptyset$, then
$
 \len_i\ge3$.
\end{claim}
\begin{proof}
By Proposition \ref{prop:attachments}, we have
 $\ell_i\geq 2$.
Suppose by contradiction that $\ell_i=2$.
If there are distinct vertices $a\in N_H(u_i)$ and $b\in N_H(u_{i+1})$,
a path from $a$ to $b$ in $H$, together with $u_ia$ and $bu_{i+1}$,
is an exterior $u_i$--$u_{i+1}$ path of length at least three,
contradicting the maximality of $C$. Since both neighborhoods
are nonempty, it follows that $N_H(u_i)=N_H(u_{i+1})=\{y\}$
for some $y\in V(H)$.
Moreover, as $S_i\cap S_{i+1}=\varnothing$,
we must have $y \notin R$.
Consequently, each family $\mathcal F_j$ contains at most one path ending
in $\{u_i,u_{i+1}\}$.
In every path in $\mathcal{F}$ ending at $u_{i+1}$,
replace its last edge $yu_{i+1}$ by $yu_i$, obtaining a new collection $\mathcal F'$.
Its set of ends on $C$ is $F\setminus\{u_{i+1}\}$,
contradicting the minimality of $|F|$.
\end{proof}

Use the weights $\lambda$ over $R$ defined as follows
\[
 \lambda_0=\lambda_3=\frac32,\qquad \lambda_1=\lambda_2=1.
\]
As usual, for any $A \subseteq \{0,1,2,3\}$, we set
 $\lambda(A)=\sum_{j\in A}\lambda_j$.

\begin{claim}\label{cl:lilambdas}
    For every $i$, $\ell_i \geq \frac{\lambda(S_i)+\lambda(S_{i+1})}{2}$.
\end{claim}
\begin{proof}
Write $U =S_i\cup S_{i+1} $ and $I =S_i\cap S_{i+1}$.
    If $S_i \cap S_{i+1}= \emptyset$ then
    $\lambda(S_i)+\lambda(S_{i+1}) \leq 5$.
    As $\ell_i \geq 3$ by Claim \ref{cl:shortbad}, the claim follows.
    Let
$
d=\max\{|a-b|:a\in I,\ b\in U\}.
$
By Claim~\ref{cl:overlap}, $\len_i\ge2+d$.

Note that
\[
\lambda(S_i)+\lambda(S_{i+1})=\lambda(U)+\lambda(I) \le 2 \lambda(U) \le 10.
\]

As $\lambda(U)+ \lambda(I) \le 10$,
we have
$\lambda(U)+ \lambda(I) \le 2 \ell_i$ when $d=3$. So assume $d<3$.
Moreover,
if $d=0$ then 
$|U|=|I|=1$. So $\lambda(U)+ \lambda(I) \le 3 \le 2\ell_i$.
Hence $1\le d\le 2$, and $\ell_i\ge 3$.
If $0\in I$ and $d=1$, then $U\subseteq\{0,1\}$, so
$\lambda(U)+\lambda(I)\le5\le2\ell_i$.
If $0\in I$ and $d=2$, then $U\subseteq\{0,1,2\}$, so
$\lambda(U)+\lambda(I)\le7\le2\ell_i$.
The case $3\in I$ is symmetric.
If $I\in\{\{1\},\{2\}\}$, then
$\lambda(U)+\lambda(I)\le5+1=6\le2\ell_i$.
It remains that $I=\{1,2\}$.
If $d=1$, then $U=I$ and the total weight is four.
If $d=2$, then
$\lambda(U)+\lambda(I)\le5+2=7\le2\ell_i$.
\end{proof}

Since $|\mathcal F_j| \geq k-3$ for each $j \in \{0,1,2,3\}$,
$ \sum_{i}\lambda(S_i)
\ge
(\frac{3}{2}+1+1 + \frac{3}{2})(k-3)
=
5(k-3).
$
By Claim \ref{cl:lilambdas},
$
|C|=\sum_i\len_i\ge\sum_i\lambda(S_i)\ge5(k-3).
$
\end{proof}

\section{Proof of Lemma \ref{lem:nop4}}

We use the following consequence of Proposition~\ref{prop:routing}.
\begin{proposition}\label{prop:Lgeq4k-1-2cap+r-2t}
Let $C$ be a longest cycle in a $k$-connected graph $G$, where
$k\ge3$, and let $R$ be a nontrivial path in $G-C$ with ends
$p,r$ such that $d_{G-C}(p)=d_{G-C}(r)=1$.
Let $F=N_C(p) \cup N_C(r)$
and $J=N_C(p) \cap N_C(r)$. 
Let $t$ be the number of $F$-segments
with at least one end in $J$. 
Then
$|C| \geq 4(k-1) -2|J|+|R|t$.
\end{proposition}
\begin{proof}
   Since $d_{G-C}(p)=d_{G-C}(r)=1$ and $\delta(G)\ge k$,
we have $d_C(p),d_C(r)\ge k-1$, so $|F|\ge2$.
    By Proposition \ref{prop:basica},
    every $F$-segment has length at least $2$. 
    Also, by Proposition \ref{prop:routing}, if a $F$-segment has at least one end in $N_C(p) \cap N_C(r)$, then it has length at least $2+|R|$.
    Hence, as
    $|F|+|J|=d_C(p)+ d_C(r) \geq 2(k-1)$, we have $|C| \geq 2|F|+|R|t \geq 4(k-1) -2|J|+|R|t$.
\end{proof}

\begin{lemma}\label{lem:nop4}
Let $C$ and $D$ be two longest cycles in a $k$-connected graph, with $k\ge2$.
If neither $G-C$ nor $G-D$ contains a path on four vertices, then
$|V(C)\cap V(D)|\ge k$.
\end{lemma}
\begin{proof}
If $C$ is Hamiltonian the conclusion follows, so let us assume that $C$ is not Hamiltonian.
Put $X=V(C)\cap V(D)$, $x=|X|$, and $L=|C|=|D|$.
As recalled in the introduction, $x\ge2$. If $k=2$, the conclusion
follows, so assume that $k\ge3$.
If every component of $G-C$ is a vertex then $L \leq 2x$ and by Proposition \ref{prop:attachments}, we have $L \geq 2k$, so $x \geq k$ as we want.
Hence assume that there exists
a component of $G-C$ with at least 2 vertices.

First suppose that neither $G-C$ nor $G-D$ contains a path on 3 vertices.
Each component of $D-X$ has at most two vertices, and hence $L\le3x$.
If $D-X$ has no edges, then $L\le2x$ and
Proposition~\ref{prop:attachments} gives $x\ge k$.
Hence assume that $D-X$ contains an edge, and choose such an
edge $pr$. Its component in $G-C$ is exactly this edge.
Write $F=N_C(p)\cup N_C(r)$, $J=N_C(p)\cap N_C(r)$
and
let $t$ be the number of $F$-segments
with at least one end in $J$. 
By Proposition \ref{prop:Lgeq4k-1-2cap+r-2t}, with $R=pr$, we have
\begin{equation}\label{eq:edgecount}
L\ge 4(k-1)-2|J|+t.
\end{equation}

Every vertex of $J$ belongs to $X$. Otherwise, replacing the edge $pr$
of $D$ by $p u r$, for $u\in J \setminus V(D)$, would increase the length of $D$.
Note also that $t \geq |J|$.
If $J=\varnothing$, then $t=0$ and
$3x\ge L\ge4(k-1)>3(k-1)$. Hence $x\ge k$ by integrality.
If $J=F$, then $F\subseteq X$.
By Proposition~\ref{prop:attachments},
$|F|  \geq k$,
so $x \geq k$.

It remains that $\varnothing\ne J\ne F$.
Let $h\ge1$ be the number of maximal sequences of consecutive
vertices of $F$ that belong to $J$. A sequence with $r$ vertices
is incident with $r+1$ $F$-segments, so $t=|J|+h$.
For each such sequence, let $u$ be its last vertex and let $v$
be the next vertex of $F$ in the orientation of $C$.
Then $u\in J$, $v\in F\setminus J$, and $|C[u,v]|\ge3$.
The set $V(C[u,v])\setminus\{u\}$ contains a vertex of $X$;
otherwise $G-D$ contains a path on three vertices.
These sets are pairwise disjoint and contain no vertex of $J$.
Thus $x\ge|J|+h=t$.
Since $|J|\le t-1$, \eqref{eq:edgecount} now gives
$L \ge 4k-4-2|J|+t \ge4k-2-t \ge4k-2-x $.
Together with $L\le3x$, this implies $4x\ge4k-2$ and hence $x\ge k$ by integrality, as we want.

It remains to consider the case in which at least one of $G-C$
and $G-D$ contains a path on three vertices
(but not four vertices).
Every subpath of $D-X$ has at most three vertices, so  $L\le4x$.
Interchanging $C,D$ if needed,
choose a component $H$ of $G-C$ containing a path $pqr$.
A connected graph containing a three-vertex path but no four-vertex path
is either a star or a triangle.
First suppose that $H$ is a star.
By Proposition \ref{prop:Lgeq4k-1-2cap+r-2t}, we have
\begin{equation}\label{eq:starcount}
L\ge 4(k-1)-2|J|+2t,
\end{equation}
where $J=N_C(p)\cap N_C(r)$,
$F=N_C(p)\cup N_C(r)$
and $t$ is the number of $F$-segments with at least one end in $J$.

If $t>|J|$, then, by \eqref{eq:starcount}, we have
$4x \geq 4k-2$ and the proof follows by integrality. Hence, assume that $t=|J|$, so by \eqref{eq:starcount}
we have $x \geq k-1$. Suppose by contradiction that $x<k$, so $x=k-1$.
Since $4(k-1)\le L\le4x$, we have $L=4(k-1)$.
If $J=\varnothing$, then $N_C(p)$ and $N_C(r)$ are nonempty
and disjoint, so there is an $F$-segment with one end in each
set. This segment has length at least four, whereas all other
$F$-segments have length at least two. Hence
$
 L\ge2|F|+2=2d_C(p)+2d_C(r)+2\ge4k-2,
$
a contradiction. Thus $J\ne\varnothing$.
Since $t=|J|$, we must have
$J=F$. Every $F$-segment then has length at least four.
As $|F|\ge k-1$ and $L=4(k-1)$, it follows that $|F|=k-1$
and every $F$-segment has length exactly four.
Proposition~\ref{prop:attachments} gives $|N_C(H)|\ge k$, so
there exists $w\in N_C(H)\setminus F$ in the interior of an
$F$-segment $C[u,v]$.
As $w$ is adjacent to a vertex of $H$, Proposition~\ref{prop:basica}
gives $|C[u,w]|,|C[w,v]|\ge3$, contradicting $|C[u,v]|=4$.

Now suppose $H$ is a triangle, with vertices $p_0,p_1,p_2$.
For any $u  \in N_C(H)$, let
$S(u)=\{j:up_j\in E(G)\}$.
Set $F=N_C(H)$, $z=|F|$, and
$I=\sum_{u\in F}|S(u)|$.
Proposition~\ref{prop:attachments} and minimum degree give
\begin{equation}\label{eq:triangleincidences}
 z\ge k,\qquad I \ge3(k-2).
\end{equation}
We say that a vertex $u \in F$ is \emph{monochromatic} if $|S(u)|=1$ and a $F$-segment is \emph{monochromatic} if its ends $u$ and $v$ are monochromatic and $S(u)=S(v)$. 
Let $s$ be the number of monochromatic vertices and $h$ the
number of monochromatic segments. Note that $h\le s$.
Also, according to Proposition~\ref{prop:basica}, every monochromatic
segment has length at least two and every other segment has length
at least four.  We have
\begin{equation}\label{eq:trianglelength}
 L\ge4z-2h.
\end{equation}
If $s=0$ then $h=0$. So \eqref{eq:triangleincidences}
and \eqref{eq:trianglelength} give $4x \ge L \ge 4k$, as we want.
Since $d_C(p_j)\ge k-2>0$ for every $j\in\{0,1,2\}$, it is
not possible that all vertices of $F$ are monochromatic with
the same set $S(u)$. As $s>0$, counting the maximal sequences
of consecutive monochromatic vertices with equal sets gives
$0\le h\le s-1$.
Also $I\le 3(z-s)+s=3z-2s$. Consequently,
\begin{equation}\label{eq:trianglechain}
 4x\ge L\ge4z-2h\ge4z-2s+2\ge z+I+2\ge4k-4.
\end{equation}
Suppose by contradiction that $x<k$, so $x=k-1$.
Equality is forced throughout in \eqref{eq:trianglechain}. In particular,
$
 z=k, I=3(k-2),  I=3z-2s,  h=s-1,
$
so $s=3$ and $h=2$.
Hence, there exist three consecutive vertices $u,v,w$ in $F$
with $S(u)=S(v)=S(w)=\{i\}$ for some $i\in\{0,1,2\}$.
Let $j\ne i$. As $d_C(p_j)\ge k-2$, the pigeonhole principle
implies that $p_j$ is adjacent to some vertex in $\{u,v,w\}$,
a contradiction.
\end{proof}

\paragraph{Note added before submission.}
During the final preparation of this manuscript, we became aware of
earlier work by Cui Wang on Smith's conjecture, recorded in
Scopus.\footnote{\url{https://www.scopus.com/pages/publications/85219094961}.}
We have not been able to obtain the full text of that article.
We have, however, consulted Chapter~3 of Wang's master's
thesis~\cite{WangThesis}, where Theorem~3.1 states that Smith's
conjecture holds when $n\le9k-36$.
We have not been able to verify all the steps of the proof
presented there, as a preliminary reading suggests inconsistencies
in some of the choices made and possible errors in some of the
arguments.
The present work was developed
independently and provides a shorter, self-contained proof, while
improving the sufficient condition to $n\le9k-30$.


\begin{thebibliography}{99}
\bibitem{BM}
J.~A. Bondy and U.~S.~R. Murty,
\emph{Graph Theory}, Graduate Texts in Mathematics 244,
Springer, 2008.

\bibitem{DChen}
D.~M. Chen,
An improved bound for Smith's longest cycles conjecture via a forbidden
subdivision,
preprint, arXiv:2608.30353, 2026.
\url{https://arxiv.org/abs/2608.30353}.

\bibitem{CFG}
G.~Chen, R.~J. Faudree, and R.~J. Gould,
Intersections of longest cycles in $k$-connected graphs,
\emph{Journal of Combinatorial Theory, Series B} 72 (1998), 143--149.

\bibitem{CCP}
E.-K. Cho, I.~Choi, and B.~Park,
Improvements on Hippchen's conjecture,
\emph{Discrete Mathematics} 345 (2022), 113029.

\bibitem{GLSTY}
C.~Groenland, S.~Longbrake, R.~Steiner, J.~Turcotte, and L.~Yepremyan,
Longest cycles in vertex-transitive and highly connected graphs,
\emph{Bulletin of the London Mathematical Society} 57 (2025), 2975--2990.
\url{https://doi.org/10.1112/blms.70134}.

\bibitem{Gro84}
M.~Gr\"otschel,
On intersections of longest cycles,
in B.~Bollob\'as (ed.), \emph{Graph Theory and Combinatorics},
Academic Press, London, 1984, pp.~171--189.

\bibitem{GN84}
M.~Gr\"otschel and G.~L. Nemhauser,
A polynomial algorithm for the max-cut problem on graphs without long
odd cycles,
\emph{Mathematical Programming} 29 (1984), 28--40.

\bibitem{G2021}
J.~Guti\'errez,
On the intersection of two longest paths in $k$-connected graphs,
\emph{Pro Mathematica} 31 (62) (2021), 11--23.

\bibitem{GV}
J.~Guti\'errez and C.~Valqui,
On two conjectures about the intersection of longest paths and cycles,
\emph{Discrete Mathematics} 347 (2024), 114148.
\url{https://doi.org/10.1016/j.disc.2024.114148}.

\bibitem{Hippchen}
T.~Hippchen,
\emph{Intersections of longest paths and cycles},
Master's thesis, Georgia State University, 2008.

\bibitem{MNZ}
J.~Ma, B.~Ning, and Z.~Zhao,
Longest cycles intersect linearly in highly connected graphs,
preprint, arXiv:2609.20724, 2026.
\url{https://arxiv.org/abs/2609.20724}.

\bibitem{MZ}
J.~Ma and Z.~Zhao,
Intersections of longest cycles in vertex-transitive and highly connected
graphs,
preprint, arXiv:2508.17438, 2025.
\url{https://arxiv.org/abs/2508.17438}.

\bibitem{SZZ}
A.~Shabbir, C.~T. Zamfirescu, and T.~I. Zamfirescu,
Intersecting longest paths and longest cycles: A survey,
\emph{Electronic Journal of Graph Theory and Applications} 1 (2013),
56--76.

\bibitem{ST95}
I.~A. Stewart and B.~Thompson,
On the intersections of longest cycles in a graph,
\emph{Experimental Mathematics} 4 (1995), 41--48.

\bibitem{Wu2026}
H.~Wu,
Intersection of cycles and paths in $k$-connected graphs,
\emph{Discrete Applied Mathematics} 378 (2026), 226--233.
\url{https://doi.org/10.1016/j.dam.2025.07.013}.

\bibitem{WangThesis}
C.~Wang,
\emph{Research on Certain Problems Concerning Cycle Length
Distribution in $k$-Connected Graphs},
Master's thesis, Jiangsu University, 2025 (in Chinese).

\end{thebibliography}
\end{document}